\documentclass[a4paper, 
11pt
]{article}

\usepackage{subfig}
\usepackage{pgfplots}
\usepackage{pgfplotstable}
\usepackage{mathtools}
\usepackage{multicol}
\usepackage{comment}
\usepackage{booktabs}
\pgfplotsset{compat=1.5}
\usepackage{amssymb}
\usepackage{url}
\usepackage{bm}
\usepackage{symbols}

\usepackage{pdfsync}
\usepackage{float}
\usepackage{tabularx}
\usepackage{enumerate}
\usepackage{array}
\usepackage{xspace}

\usepackage{authblk}
\usepackage{hyperref}

\newcommand{\TODO}[1]{{\color{black}#1}}

\usepgfplotslibrary{patchplots} \usetikzlibrary{shapes.geometric,
  arrows} \pgfplotsset{ compat=1.6, standard/.style={ axis x
    line=middle, axis y line=middle, enlarge x limits=0.15, enlarge y
    limits=0.15, every axis x label/.style={at={(current axis.right of
        origin)},anchor=north west}, every axis y
    label/.style={at={(current axis.above origin)},anchor=north east}
  } }

\usepackage{booktabs}
\pgfplotstabletypeset[
	every head row/.style={
		before row=\toprule,after row=\midrule},
	every last row/.style={
		after row=\bottomrule},
]

\newcommand{\beginpipe}[2]
{
  \draw (#1    , 0  ) --
        (#2-0.5, 0  ) --
		(#2,     0.6) --
		(#2-0.5, 1.2  ) --
		(#1,     1.2  ) --
  		(#1,     0  )

}

\newcommand{\pipe}[3]
{
	\filldraw[color=black, fill=#3] (#1-0.5, 0) --
        (#2-0.5, 0) --
		(#2,     .6) --
		(#2-0.5, 1.2) --
		(#1-0.5, 1.2) --
		(#1,     .6) --
  		(#1-0.5, 0)

}

\newcommand{\epipe}[2]
{
  \draw (#1-0.5, 0  ) --
        (#2,     0  ) --
		(#2,     1.2  ) --
		(#1-0.5, 1.2  ) --
		(#1,     0.6) --
  		(#1-0.5, 0  )

}

\begin{document} 

\title{Dimension reduction in heterogeneous parametric spaces with application
  to naval engineering shape design problems}

\author[]{Marco~Tezzele\footnote{marco.tezzele@sissa.it}}
\author[]{Filippo~Salmoiraghi\footnote{fsalmoir@sissa.it}}
\author[]{Andrea~Mola\footnote{andrea.mola@sissa.it}}
\author[]{Gianluigi~Rozza\footnote{gianluigi.rozza@sissa.it}}

\affil[]{Mathematics Area, mathLab, SISSA, International School of Advanced Studies, via Bonomea 265, I-34136 Trieste, Italy}

\maketitle

\begin{abstract}
We present the results of the first application in the naval
architecture field of a methodology based on active subspaces properties
for parameter space reduction. The physical problem considered is the
one of the simulation of the hydrodynamic flow past the hull of a ship advancing in calm
water. Such problem is extremely relevant at the preliminary stages of the ship
design, when several flow simulations are typically carried out by the engineers
to assess the dependence of the hull total resistance on the geometrical parameters
of the hull, and others related with flows and hull properties. Given the high number of geometric and physical parameters which might
affect the total ship drag, the main idea of this work is to employ the active subspaces
properties to identify possible lower dimensional structures in the parameter space.
Thus, a fully automated procedure has been implemented to produce
several small shape perturbations
of an original hull CAD geometry, in order to exploit the resulting
shapes and to run high fidelity flow
simulations with different structural and physical parameters as well,
and then collect data for the active subspaces analysis. The free form deformation
procedure used to morph the hull shapes, the high fidelity solver based on potential
flow theory with fully nonlinear free surface treatment, and the active subspaces
analysis tool employed in this work have all been developed and integrated within
SISSA \emph{mathLab} as open source tools.
The contribution will also discuss several details of the implementation of such tools,
as well as the results of their application to the selected target engineering
problem.
\end{abstract}



\section{Introduction}
\label{sec:intro}

\begin{figure}[h]
\centering
\begin{tikzpicture}[node distance = 2.5cm, auto]
	\beginpipe{0}{2.5};
	\pipe{2.55}{4.1}{white};
	\pipe{4.15}{5.7}{white};
	\pipe{5.75}{7.3}{white};
	\pipe{7.35}{8.9}{white};
        \epipe{8.95}{11.3};
        
	\node[align=left] at (1.0,.6) {Problem\\settings};
	\node[align=left] at (3.2,.6) {FFD};
	\node[align=left] at (4.8,.6) {BEM};
	\node[align=left] at (6.3,.6) {AS};
	\node[align=center] at (7.9,.6) {RSM};
	\node[align=center] at (10.,.6) {Fast CFD\\evaluations};
\end{tikzpicture}
\caption{Scheme of the structure of the pipeline proposed. The
  geometrical deformation is performed via free from deformation (FFD),
  then a Boundary Elements Method (BEM) solver computes the wave resistance, the active subspaces
  (AS) are detected and finally a response surface method (RSM) is employed.}
\label{fig:scheme}
\end{figure}
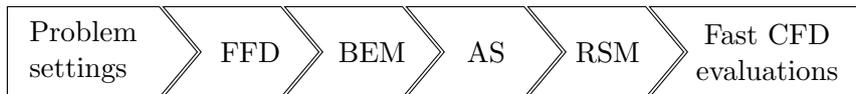

Nowadays engineering simulations present a wide range of different
parameters. When it comes to find an optimal solution with respect to
the physical constraints it is easy to be affected by the curse of
dimensionality, when the number of parameters makes the simulation
unfeasible. This problem arises quite easily even with a small
parameter space dimension (depending on the simulation, even ten
parameters could take months to be optimized). In this framework
reducing the dimension of this space becomes crucial and a priority. 
To tackle it we focus on the active subspaces property (see~\cite{constantine2015active}) to carry out a technique applied on a naval
engineering problem, that is the computation of the total wave
resistance of a hull advancing in calm water. In the framework of simulation-based
design and shape optimization we cite, among
others,~\cite{diez2015design,tahara2012cfd,volpi2014development,dambrine2016theoretical}. The
computational pipeline we are
going to present is composed first by a geometrical parametrization and deformation of the hull
through free form deformation (see~\cite{sederberg1986free}). Then the use of a high fidelity solver
based on Boundary Elements Method (BEM) to get the wave resistance with
respect to the geometrical parameters. We consider also a structural parameter --- the
initial displacement of the hull --- and a physical one --- the
velocity of the hull ---. Subsequently active subspaces are
identified thanks to the data collected from the high fidelity solver, and
finally a proper reduced response surface is constructed. The result allows
the final user to have an estimate of the wave resistance under a
certain threshold and within a time of one second with respect to the hours
needed for a single classic full simulation. Moreover during the process is
possible to identify the most important parameters and have insights
on how they influence the output of interest. Figure~\ref{fig:scheme}
summarizes the proposed computational pipeline.

The content of this contribution is organized as follows. 
Section \ref{sec:potential} introduces the ship resistance prediction problem, its
dependence on hull shape deformations, and equations of the fluid structure interaction
model used for the simulations.
In section~\ref{sec:ffd} we recall the free form
deformation technique and we show the main features of the developed
tool to manage parametric shapes. Section~\ref{sec:bem} has the purpose of introducing
some detail about the high fidelity solver implementation. 
In section~\ref{sec:active} we present the active subspaces properties and
its features, with a numerical recipe to identify them.
Then section~\ref{sec:results} shows the numerical results obtained by coupling the
three methods in sequence. Finally conclusions and perspectives are
drawn in section~\ref{sec:the_end}.

\section{A model naval problem: wave resistance estimation of a hull advancing in calm water}
\begin{figure}[htb]
\centering
\includegraphics[width=1.\textwidth]{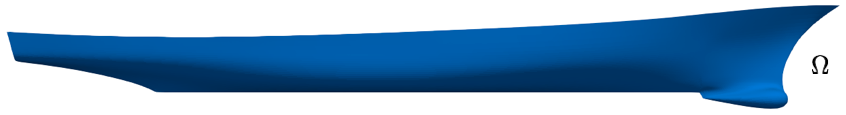}
\caption{Representation of the reference domain $\Omega$, that is the DTMB 5415 hull.}
\label{fig:original_hull}
\end{figure}

In this section we introduce the problem of the estimation of the
resistance of a hull advancing in calm water. The model hull shape considered in
this work is the DTMB 5415, which was originally conceived for the preliminary design of
a US Navy Combatant ship. Due to the wealth of experimental data available
in the literature (see for example~\cite{olivieri2001towing,stern2001international}) such shape, which includes
a sonar dome and a transom stern (see Figure~\ref{fig:original_hull}), has become a common
benchmark for naval hydrodynamics simulation tools.

Let $\Omega \subset \mathbb{R}^3$, be a domain (see Figure
\ref{fig:original_hull}) associated with our DTMB 5415 model hull.  
We call $\Omega$ the \emph{reference} domain; for practical reasons this
domain happens to correspond to the undeformed hull, even
though this assumption is not fundamental for the remainder of the
paper. We here remark that the domain considered in the fluid dynamic simulations
is in principle the volume surrounding the hull which is occupied by water,
namely $\Omega_w$. Further details about the fluid dynamic
domain will be provided in Section~\ref{sec:potential}.

Let $\mathcal{M}(\boldsymbol{x}; \boldsymbol{\mu}^{\text{GEOM}}): \mathbb{R}^3 \to \mathbb{R}^3$ be
a shape morphing that maps the reference domain $\Omega$ into the deformed
domain $\Omega(\boldsymbol{\mu}^{\text{GEOM}})$ as follows:
\begin{equation*}
\Omega(\boldsymbol{\mu}^{\text{GEOM}}) = \mathcal{M}(\Omega; \boldsymbol{\mu}^{\text{GEOM}}).
\end{equation*}
Quite naturally, the results of the fluid dynamic simulations will depend on the
specific hull shape considered, which are in turn associated to the parameters defining the
morphing $\mathcal{M}$ (which will be exensively  defined in Section~\ref{sec:ffd}).
It is worth pointing out here that the geometrical quantity having the most effect on the resistance is
the immersed volume of a hull shape, as higher volumes will generate higher drag values. \TODO{This is clearly due to
the fact that higher hull volumes will result in a higher mass of water displaced as the ship advances in the
water, and in increased surface exposed to the water friction. This consideration might lead to the naive
conclusion that since the shape optimizing the total drag is the one corresponding to zero buoyant volume,
the hull volume must be constrained in the optimization algorithm, to avoid a convergence to such trivial shape, which would
not generate a vertical force able to sustain the ship weight.
In this work, rather than constraining the hull volume through more complex hull deformation algorithms, we decided to impose the
weight (or displacement) of the hull in the fluid dynamic simulations, so that each hull would reach its hydrostatic
equilibrium position, in which the weight prescribed at the design stage is balanced by the vertical hydrodynamic force.
This solution, which of course requires a model accounting for the rigid motions of the hull into the fluid dynamic simulations,
is able to lead to design solutions which optimize the total resistance while retaining the required load capability of the ship.} 
   
For all the aforementioned reasons, along with the geometrical parameters associated with the hull morphing, the results of
our simulations are also affected by the ship displacement and cruise speed, which are instead physical
parameters determining the hydrodynamic equilibrium position and forces.
Thus, considering both the geometric morphing and variations in the physical parameters, we have a set
of $m\in \mathbb{N}$ parameters which affect the output of the fluid dynamic simulations. The parametric
domain considered is then defined as $\mathbb{D} \subset \mathbb{R}^m$, and is assumed to be a box in $\mathbb{R}^m$. 

By a practical standpoint, once a point in the parameter domain $\mathbb{D}$ is identified, the specific hull
geometry as well as the desired ship displacement and cruise speed are provided to the fluid dynamic solver, which
carries out a flow simulation to provide a resistance estimate. In this framework free form deformation has been
employed for the generation of a very large number of hull geometries obtained from the morphing of the DTMB 5415
naval combatant hull. Each geometry generated has been used to set up a high fidelity hydrodynamic simulation with
the desired ship displacement and hull speed. The output resistances for all the configurations tested have been
finally analyzed by means of active subspaces in order to reduce the parameter space.

In the next subsections, we will provide a brief description of the
unsteady fully nonlinear potential fluid dynamic model used to carry out the
high fidelity simulations. In addition, we will describe the rigid body
equations based on hull quaternions used to compute the hull linear and
angular displacements corresponding to the final hydrodynamic equilibrium
position reached at the end of each simulation. We refer the interested
reader to~\cite{molaEtAl2013,mola2016ship,MolaHeltaiDeSimone2017} for
further information on the fully nonlinear potential free surface model,
on its application to complex hull geometries, and on the treatment of
the hull rigid motions, respectively.

\subsection{Fully nonlinear potential model}
\label{sec:potential}

In the simulations we are only considering the motion of a ship advancing
at constant speed in calm water. For such reason we solve the problem
in a \emph{global}, translating reference frame $\widehat{XYZ}$,
which is moving with the constant horizontal velocity of the
boat $\Vb_\infty = (V_\infty,0,0)$. Thus, the $X$ axis
of the reference frame is aligned with $\Vb_\infty$, the $Z$ axis is
directed vertically (positive upwards), while the $Y$ axis is directed
laterally (positive port side). 

As aforementioned, the domain $\Omega_w(t)$ in which we are interested in
computing the fluid velocity and pressure is represented by the portion of
water surrounding the ship hull. The time varying shape of such
domain --- and in particular that of its boundary with the air above --- is
one of the unknowns of the fluid dynamic problem. By convention, we place the
origin of the vertical axis $Z$ in correspondence with the undisturbed free surface
level, and we start each simulation at time $t=0$ from such undisturbed
configuration. \TODO{Thus, at least in its initial configuration, the flow domain is represented by
$\Omega_w(t=0) = \mathbb{R}^3_{Z-}\backslash \Omega$, which is the boolean subtraction of the hull volume
from the lower half-space of $\mathbb{R}^3$ for which $Z\leq 0$, here indicated with $\mathbb{R}^3_{Z-}$.}

If overturning waves are not observed (which is typically the case for low cruise velocities
typical of a ship), the domain $\Omega_w(t)$ is simply connected. So under the assumptions of
irrotational flow and non viscous fluid the velocity field $\vb(\Xb,t)$
admits a representation through a scalar potential function
$\Phi(\Xb,t)$, namely
\begin{equation}
\label{eq:potential-definition}
\vb = \nablab\Phi = \nablab (\Vb_\infty\cdot\Xb + \phi ) \qquad\qquad
\forall \ \Xb\in \Omega_w(t),
\end{equation}
in which $\phi(\Xb,t)$ is the \emph{perturbation potential}.
Under the present assumptions, the equations of motion simplify to the unsteady
Bernoulli equation and to the Laplace equation for the perturbation
potential:
\begin{subequations}
  \label{eq:incompressible-euler-potential}
  \begin{alignat}{3}
    \label{eq:bernoulli}
    & \frac{\Pa\phi}{\Pa t}+\frac{1}{2}\left|\nablab\phi+\Vb_\infty\right|^2
      +\frac{p-p_a}{\rho}-\gb\cdot\Xb
       = C(t) \qquad & \text{ in } \Omega_w(t) ,\\
    \label{eq:incompressibility-potential}
    & \Delta \phi = 0 & \text{ in } \Omega_w(t) ,
  \end{alignat}
\end{subequations} 
where $C(t)$ is an arbitrary function of time, and $\gb = (0,0,-g)$,
is the gravity acceleration vector, directed along the $z$ axis.
The unknowns of such mathematical problem $\phi$ and $p$ are
uncoupled. This means that the solutuon of the Poisson problem in Eq.~\eqref{eq:incompressibility-potential}
can be obtained independently of the pressure field. Once such solution
is obtained, the pressure can be obtained through a  postprocessing step
based on Bernoulli Eq.~\eqref{eq:bernoulli}. Thus, the Laplace equation is the
governing equation of our model. Such equation is complemented by non penetration
boundary conditions on the hull surface $\Gamma^{b}(t)$ and water basin bottom boundary
$\Gamma^{bot}(t)$, and by homogeneous Neumann boundary conditions on the truncation
boundaries $\Gamma^{far}(t)$ of the numerical domain. The bottom of the basin is located
at a depth corresponding to 2 boat lenghts, while the truncation boundaries are
located approximatively at a distance from the boat of 6 boat lenghts in the longitudinal
direction $X$ and of 2 boat lengths in the lateral direction $Y$). 
On the water free surface $\Gamma^{w}(t)$, we employ the kinematic and dynamic Semi-Lagrangian fully nonlinear boundary
conditions, which respectively read
\begin{eqnarray}
\frac{\delta\eta}{\delta t} &=& \frac{\Pa\phi}{\Pa z}+\nablab\eta\cdot\left(\wb -\nablab\phi-\Vb_\infty\right)
\qquad \text{ in } \Gamma^{w}(t) ,
\label{eq:fsKinematicBeck}\\
\frac{\delta\phi}{\delta t} &=& -g\eta +
\frac{1}{2}|\nablab\phi|^2 + \nablab\phi\cdot\left(\wb -\nablab\phi-\Vb_\infty\right)
\qquad \text{ in } \Gamma^{w}(t).
\label{eq:fsDynamicBeck}
\end{eqnarray}  
The former equation expresses the fact that a material point moving on the
free surface will stay on the free surface --- here assumed to be
a single valued function $\eta(X,Y,t)$ of the horizontal
coordinates $X$ and $Y$. The latter condition
is derived from Bernoulli Eq.~\eqref{eq:bernoulli},
under the assumption of constant atmospheric pressure on the water surface.
This peculiar form of the fully nonlinear boundary
conditions was proposed by Beck \cite{beck1994}. Eq.~\eqref{eq:fsKinematicBeck}
allows for the computation of the vertical velocity of markers which move on the
water free surface with a prescribed horizontal speed $(w_X,w_Y)$.
Eq.~\eqref{eq:fsDynamicBeck} is used to obtain the velocity potential
values in correspondence with such markers. The resulting
vector $\wb=(w_X,w_Y,\frac{\delta\eta}{\delta t}) = \dot{\Xb}$ is the time
derivative of the position of the free surface markers.
In this work, such free surface markers are chosen as the free surface nodes
of the computational grid. To avoid an undesirable mesh nodes drift
along the water stream, the markers arbitrary horizontal velocity is
set to 0 along the $X$ direction. The $Y$ component of the water nodes in
contact with the ship --- which is moved according with the computed linear and
angular displacements --- is chosen so as to keep such nodes on the hull surface.
As for the remaining water nodes, the lateral velocity value is set  
to preserve mesh quality.

\subsection{Three dimensional hull rigid motions}
\begin{figure}
\includegraphics[width=1.\textwidth]{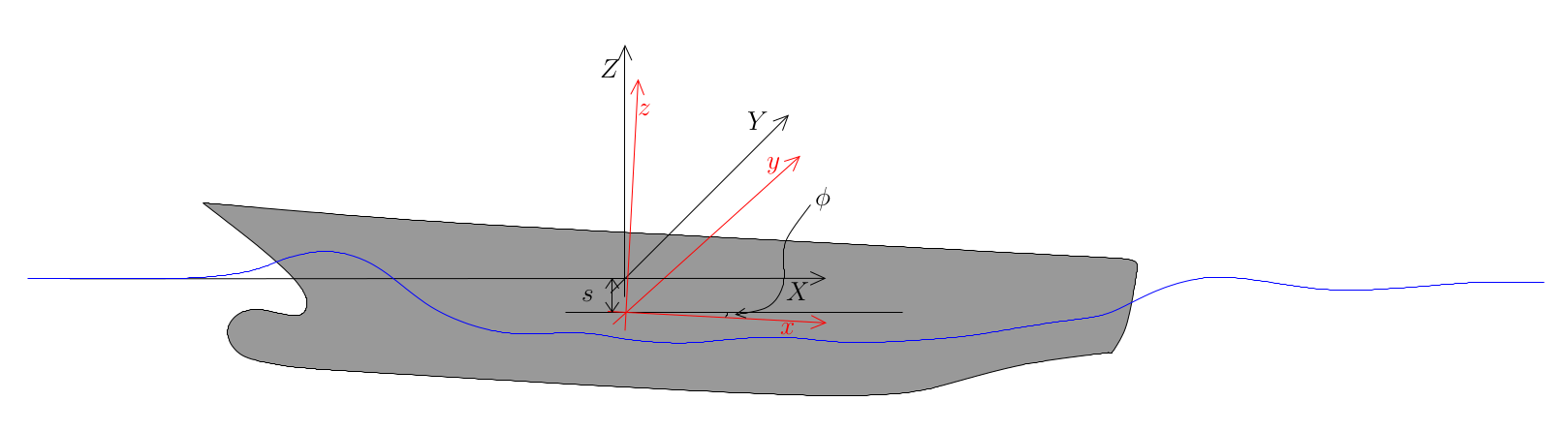}
\caption{\TODO{A sketch illutrating the hull-attached frame $\widehat{xyz}$
  in red and the global reference frame $\widehat{XYZ}$ which is
  moving with the constant horizontal velocity of the boat. The ship
  here depicted is experiencing a vertical displacement $s$ and an
  angular displacement characterized by the pitch angle $\phi$.}}
\label{fig:sink_and_pitch}
\end{figure}

The ship hull is assumed to be a rigid body. A second, \emph{hull-attached}
reference frame $\widehat{xyz}$, which follows the hull in its translations and
rotations is employed to study the ship motions. The center of such reference
frame is located at the ship center of gravity, which in the global
reference frame reads $\Xb^G(t) = X^G(t)\eb_X+Y^G(t)\eb_Y+Z^G(t)\eb_Z$,
where $\eb_X$, $\eb_Y$, $\eb_Z$ are the unit vectors along the global
system axes. See Figure~\ref{fig:sink_and_pitch} for a detailed sketch.

The rotation matrix $R(t)$ is used to convert the coordinates of
a point $\xb$ written in the hull-attached reference frame, to those in
the global frame $\Xb$, namely
\begin{equation}
\label{eqHullPointPos}
\Xb(t) = R(t)\xb+\Xb^G(t).
\end{equation}
The global frame velocity of a point having coordinates $\xb$
in the hull-attached frame is obtained as
\begin{equation}
\label{eqHullPointVel}
\Vb^{\text{hp}}(t) = \omegab(t)\times\xb+\dot{\Xb}^G(t),
\end{equation}
in which $\omegab$ is the angular velocity vector.

Eqs.~\eqref{eqHullPointPos} and \eqref{eqHullPointVel} imply that once
$\Xb^G(t)$, $R(t)$, and $\omegab(t)$ are known at time $t$, the position
and velocity of each point of the hull can be obtained. For this reason,
writing the time evolution equations for $\Xb^G(t)$, $R(t)$, and $\omegab(t)$ is
sufficient to fully determine the hull dynamics. In the next sections,
we will present such evolution equations written in the global reference frame,
as presented in \cite{Azueta2001c} and \cite{FormaggiaEtAlRowing2010}.

The evolution equation for $\Xb^G(t)$ is obtained via the linear momentum
conservation equation, which in the case of our hydrodynamics simulation
framework reads
\begin{equation}
\label{eqHullLinCons}
m_s\ddot{\Xb}^G(t) = m_s\gb + \Fb^w(t).
\end{equation}
In Eq.~\eqref{eqHullLinCons}, $m_s$ is the mass of the ship, while the
hydrodynamic force vector $\Fb^w(t)$ is in the present work obtained as
the sum of the pressure and viscous forces on the hull. 

The evolution equation for $\omegab$ is obtained writing the angular
momentum conservation, namely
\begin{equation}
\label{eqHullAngCons}
R(t) I^G R(t)^T\dot{\omegab}(t) + \omegab(t) \times R(t) I^G R(t)^T \omegab(t) = \Mb^w(t),
\end{equation}

where $I^G$ is the matrix of inertia of the ship in the hull-attached
reference frame, and hydrodynamic moment vector $\Mb^w(t)$ is the sum of the
moment about the ship center of gravity of the pressure and viscous forces
on hull, propeller and appendages.

To write an evolution equation for $R$, starting from the 
angular velocity vector $\omegab(t)=[\omega_X(t),\omega_Y(t),\omega_Z(t)]$,
we first introduce the skew symmetric tensor

\begin{equation}
\label{eqAngVelTens}
\omega(t) =
\left[\begin{array}{c c c}
0 & -\omega_Z(t) & \omega_Y(t) \\
\omega_Z(t) & 0 & -\omega_X(t) \\
-\omega_Y(t) & \omega_X(t) & 0
\end{array}\right].
\end{equation}

Note that tensor $\omega(t)$ will act on a vector $\ub\in \mathbb{R}^3$ as if
the vector product by $\omegab(t)$ were applied to $\ub$:
\begin{equation}
\omegab(t)\times\ub = \omega(t)\ub.
\end{equation}

Making use of such tensor, an evolution equation for the rotation matrix $R(t)$
reads
\begin{equation}
\label{eqRotMatEvo}
\dot{R}(t) = \omega(t) R(t),
\end{equation}

which can be advanced in time to obtain the components of $R$ and close the
equations of motions of a rigid body in three dimension. Yet, in the common
practice of rigid body simulations, direct numerical integration of
Eq.~\eqref{eqRotMatEvo} is avoided. The most important reason for this
is related to numerical drift. If we in fact keep track of the orientation
of a rigid body integrating Eq.~\eqref{eqRotMatEvo}, numerical error will
build up in the entries of $R(t)$, so that it will no longer be a rotation
matrix, losing its properties of orthogonality and of having determinant equal
to 1. Physically, the effect would be that applying $R(t)$ to a body would
cause a skewing effect.

A better way to represent the orientation of a rigid body in three dimensions
(even with large rotations) is represented by the use of \emph{unit quaternions}
(see the work of \cite{shoemake1985} for details). For our purposes, quaternions
can be considered as a particular type of four element vector, normalized to
unit length. If we indicate the
quaternion $\qb = s + v_X\eb_X + v_Y\eb_Y + v_Z\eb_Z $ as $\left[s,\vb\right]$,
the internal product of two quaternions $\qb_1$ and $\qb_2$ is defined as
\begin{equation}
\label{eqQuatProd}
\qb_1\qb_2 = \left[s_1,\vb_1\right]\left[s_2,\vb_2\right] =
\left[s_1s_2-\vb_1\cdot\vb_2\, , \, s_1\vb2+s_2\vb1+\vb_1\times\vb_2\right].
\end{equation}

The norm of a quaternion $\qb$ is defined as
$||\qb|| = \sqrt{s^2+v_X^2+v_Y^2+v_Z^2}$. Unit quaternions can be used
to represent rotations in a three dimensional space. In fact, given a
quaternion $\qb(t): ||\qb(t)||=1 \quad\forall t$, we can obtain the corresponding
rotation matrix as
\begin{equation}
\label{eqQuaterToRotMat}
R =
\left[\begin{array}{c c c}
1-2v_Y^2-2v_Z^2 &  2v_Xv_Y-2s v_Z & 2v_Xv_Z+2s v_Y \\
2v_Xv_Y+2s v_Z & 1-2v_Y^2-2v_Z^2 & 2v_Yv_Z-2s v_X \\
2v_Xv_Z-2s v_Y & 2v_Yv_Z+2s v_X & 1-2v_Y^2-2v_Z^2
\end{array}\right],
\end{equation} 
in which to lighten the notation we omitted the time dependence of both $R(t)$ and
the components of $\qb(t)$.

Finally, the equation needed to describe the time evolution for the hull
quaternion $\qb(t)$ is 
\begin{equation}
\label{eqQuaterEvo}
\dot{\qb}(t) = \dsfrac{1}{2}\omegab_q(t)\qb(t),
\end{equation}
where $\omegab_q(t) = \left[0,\omegab(t)\right]$ is the quaternion associated with
the angular velocity vector $\omegab(t)$. As quaternions only have four entries,
there only is one extra variable used to describe the three degrees freedoms of a
three dimensional rotation. A rotation matrix instead employs nine parameters for
the same three degrees of freedom; thus, the quaternions present far less
redundancy than rotation matrices. Consequently, quaternions experience far less
numerical drift than rotation matrices. The only possible source of drift in a
quaternion occurs when the quaternion has lost its unit
magnitude. This can be easily corrected by periodically renormalizing the quaternion
to unit length~\cite{shoemake1985}.

\section{Shape morphing based on Free Form Deformation}
\label{sec:ffd}
As already mentioned, we are interested in problems characterized by both physical and geometrical parameters.
In such framework, the Free Form Deformation (FFD) approach is adopted to implement the hull deformations corresponding
to each geometrical parameter set considered. 

A very detailed description of FFD is beyond the scope of the present work. 
In the following we will give only a brief overview. For a
further insight see~\cite{sederberg1986free} for the original formulation and~\cite{lombardi2012numerical,rozza2013free,salmoiraghi2016isogeometric,forti2014efficient,salmoiraghi2017} for more recent works.

We decided to adopt free form deformation among other possibilities (including, for instance, Radial Basis Functions or Inverse Distance Weighting) because it allows to have global deformations with a few parameters. 
For the complexity of the problem at hand, by trying to reduce the number of parameters starting from hundreds of them can be infeasible for the number of Monte Carlo simulations required.
One of the possible drawbacks of FFD is generally that the parameters do not have a specific geometric meaning, like, for instance, a prescribed length or distance. 
In the case of application to \TODO{active subspaces (AS)} this is not a problem since AS itself identifies new parameters, obtained by combination of the original ones, meaningless from the geometric and physical point of view.

\begin{figure}
\centering
\includegraphics[width=0.65\textwidth]{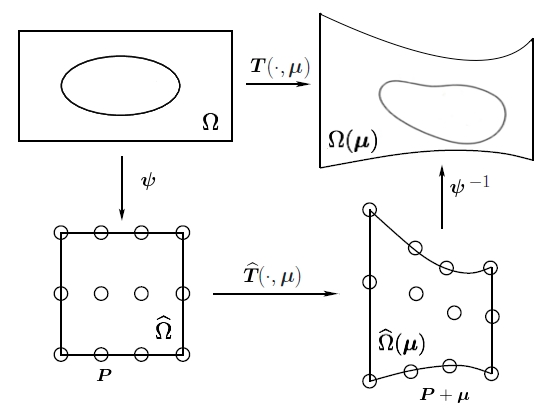}
\caption{Sketch of the FFD map construction. For ease of readability
  we dropped the superscript from $\boldsymbol{\mu}^{\text{GEOM}}$.}
\label{fig:FFD sketch}
\end{figure}

FFD consists basically in three different step, as shown in Figure~\ref{fig:FFD sketch}:
\begin{itemize}
\item Mapping the physical domain $\Omega$ to the reference one $\widehat{\Omega}$ with the map~$\boldsymbol{\psi}$.
\item Moving some control points $\boldsymbol{P}$ to deform the
  lattice with $\widehat{T}$. The movement of the control points is given
  by the weights of FFD, and represent our geometrical parameters $\boldsymbol{\mu}^{\text{GEOM}}$.
\item  Mapping back to the physical domain $\Omega(\boldsymbol{\mu})$ with the map $\boldsymbol{\psi}^{-1}$.
\end{itemize}	
	  
\begin{figure}
\centering
\includegraphics[width=0.8\textwidth]{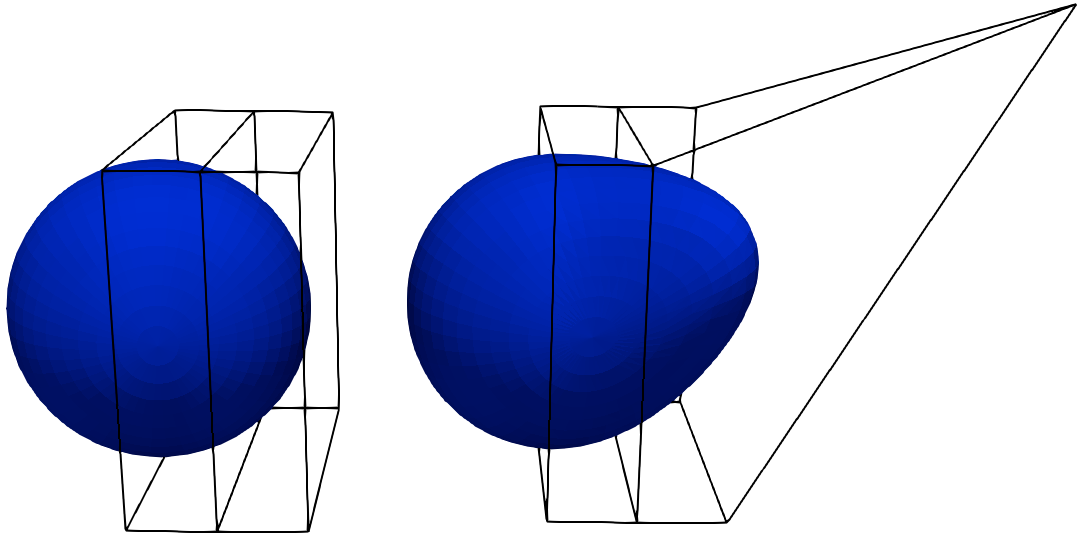} 
\caption{FFD morphing on a simple geometry: the sphere case. Here we
  move only one FFD control point in the lattice.}
\label{fig:FFD on sphere}
\end{figure}

So FFD map $T$ is the composition of the three maps, i.e.
\begin{equation}
T(\cdot, \boldsymbol{\mu}^{\text{GEOM}}) = (\boldsymbol{\psi}^{-1} \circ \widehat{T} \circ \boldsymbol{\psi})
(\cdot, \boldsymbol{\mu}^{\text{GEOM}}) .
\end{equation}
In particular, in the three dimensional case, for every point $\boldsymbol{X} \in \Omega$ inside the
FFD box, its position changes according to
\begin{equation}
T(\boldsymbol{X}, \boldsymbol{\mu}^{\text{GEOM}}) = \boldsymbol{\psi}^{-1} \left( \sum_{l=0} ^L \sum_{m=0}
  ^M \sum_{n=0} ^N
  b_{lmn}(\boldsymbol{\psi}(\boldsymbol{X}))
  \boldsymbol{P}_{lmn}^0 \left(\boldsymbol{\mu}^{\text{GEOM}}_{lmn}\right) \right) ,
\end{equation}
where $b_{lmn}$ are Bernstein polynomials of degree $l$, $m$,
$n$ in each direction, respectively, \TODO{$\boldsymbol{P}_{lmn}^0
\left(\boldsymbol{\mu}^{\text{GEOM}}_{lmn}\right) =
\boldsymbol{P}_{lmn} + \boldsymbol{\mu}^{\text{GEOM}}_{lmn}$, and
$\boldsymbol{P}_{lmn}$ represents the coordinates of the control
point identified by the three indices $l$, $m$, $n$ in the lattice of
control points. We also explicit the $\hat{T}$ mapping as follows
\begin{equation*}
\hat{T}(\boldsymbol{Y}, \boldsymbol{\mu}^{\text{GEOM}}) := \sum_{l=0} ^L \sum_{m=0}
  ^M \sum_{n=0} ^N
  b_{lmn}(\boldsymbol{Y})
  \boldsymbol{P}_{lmn}^0
  \left(\boldsymbol{\mu}^{\text{GEOM}}_{lmn}\right) \qquad \forall \,
  \boldsymbol{Y} \in \widehat{\Omega}.
\end{equation*}
}
In Figure~\ref{fig:FFD on sphere}, we show, for example, the application of the FFD
morphing on a very simple geometry, that is a sphere.

We implemented this algorithm in a \emph{python} package called \emph{PyGeM}~\cite{pygem} in order to deal with the major industrial CAD file
formats. It handles \emph{iges}, \emph{stl}, \emph{step}, \emph{vtk},
\emph{unv}, \emph{keyword}, and \emph{openfoam} files. It
extracts the coordinates of the control points, deforms them according
to the inputs given by the user and writes a new file with the morphed CAD.
We improve the traditional version of the algorithm by allowing a rotation of the FFD lattice in order to give more flexibility to the tool. In general with our package it
is possible to have a generic bounding box (not only a cube) as long
as the $\boldsymbol{\psi}$ map is affine. 

\begin{figure}[h!]
\centering
\subfloat[Front view]{\includegraphics[width=0.18\textwidth]{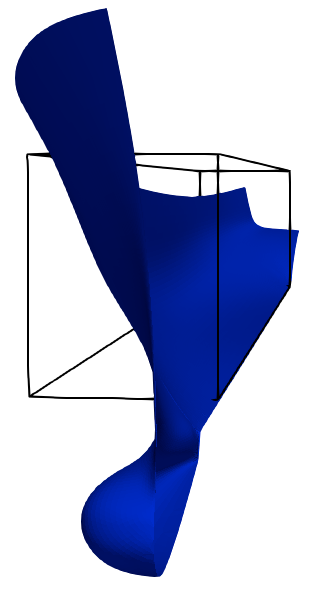}\label{subfig:ffd_hull_1}}\hspace{.2cm}
\subfloat[Back view]{\includegraphics[width=0.79\textwidth]{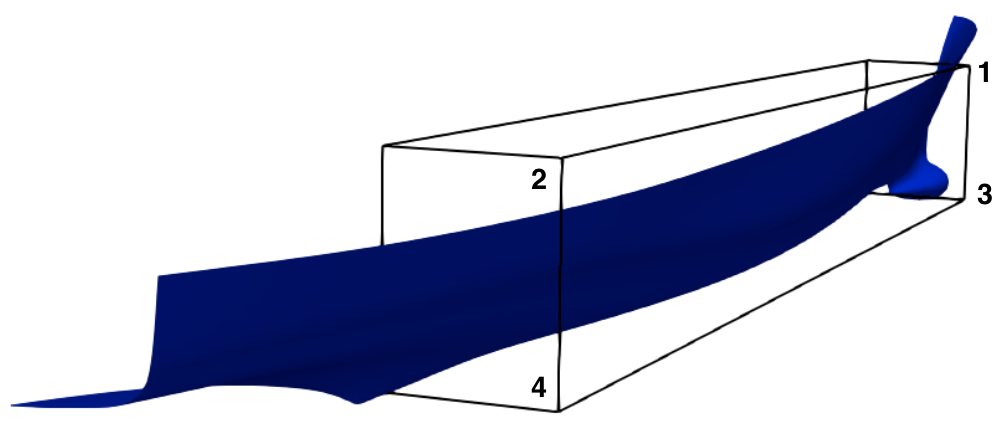}\label{subfig:ffd_hull_2}}
\caption{Plot \protect\subref{subfig:ffd_hull_1} shows the FFD lattice over
  one side wall of the hull from the front, while plot \protect\subref{subfig:ffd_hull_2}
  depicts the hull and the lattice from the back together with the
  numbers that identify the FFD points.}
\label{fig:ffd_lattice}
\end{figure}

In order to exemplify the equations above to our case, let us consider
Figure~\ref{fig:ffd_lattice}, where the control points we are going to
move are marked with numbers. As geometrical parameters we select six components of these four control
points of the FFD lattice over one side wall of the hull. Then we apply
the same deformation to the other side. This because one of our
hypothesis is the symmetry of the deformed hull. In particular
table~\ref{tab:design_space} summarizes the set of design variables,
the associated FFD-node coordinate modified ($y$ is the span of the
hull, $x$ its length and $z$ its depth) and the lower and upper bound
of the modification. There are also two more parameters that do not
affect the geometry, and are related to the physics of the problem,
that is the displacement and the velocity of the hull. From now on we denote
with $\boldsymbol{\mu} := \{ \mu_i \}_{i \in [1, \dots, 8]}$ the
column vector of the parameters, where $\mu_i$ are defined in
table~\ref{tab:design_space}. To denote only the parameters affecting
the geometrical deformation we use $\boldsymbol{\mu}^{\text{GEOM}} := \{ \mu_i \}_{i \in [1, \dots, 6]}$. For sake of clarity we underline that
the undeformed original domain is obtained setting all the geometrical
parameters to 0. All the upper and lower bounds are chosen in order to satisfy physical constraints. 

\begin{table}[h]
\centering
\caption{Design space (FFD lattice $2 \times 2 \times 2$) with eight
  design parameters. Six geometrical parameters chosen among the FFD control
points, one structural parameter that is the initial displacement of the
hull and one physical parameter given by the velocity.}
\label{tab:design_space}
\begin{tabular}{cccc}
\toprule
Parameter & Nature & Lower bound & Upper bound \\ \midrule
$\mu_1$ & FFD Point 1 $y$ & -0.2 & 0.3 \\
$\mu_2$ & FFD Point 2 $y$  & -0.2 & 0.3 \\
$\mu_3$ & FFD Point 3 $y$  & -0.2 & 0.3 \\
$\mu_4$ & FFD Point 4 $y$  & -0.2 & 0.3 \\ 
$\mu_5$ & FFD Point 3 $z$ & -0.2 & 0.5 \\
$\mu_6$ & FFD Point 4 $z$  & -0.2 & 0.5 \\ \midrule
$\mu_7$ & weight (kg) & 500 & 800 \\ \midrule
$\mu_8$ & velocity (m/s) & 1.87 & 2.70 \\ \bottomrule
\end{tabular}
\end{table}

To create the dataset with all the deformation, we started from the
original \emph{iges} file of the hull and deformed it with the
\emph{PyGeM} package. The deformations are generated randomly with an
uniform distribution.

\section{Implementation of high fidelity potential solver based on
         the Boundary Element Method}
\label{sec:bem}        

The boundary value problem described in Section~\ref{sec:potential} is governed by the linear Laplace operator. Yet, it
is nonlinear due to the presence of the boundary conditions in
Eqs.~\eqref{eq:fsKinematicBeck} and \eqref{eq:fsDynamicBeck}. Further
sources of nonlinearity are given by continuous change of the domain shape over time
and by the arbitrary shape of the ship hull. Thus, for each time instant, we will
look for the correct values of the unknown potential and node displacement fields by
solving a specific nonlinear problem resulting from the spatial and time discretization
of the original boundary value problem. The spatial discretization of the Laplace
problem is based upon a boundary integral formulation described in \cite{giulianiEtAl2015}.
In this framework, the domain boundary is subdivided into quadrilateral cells, on
which bi-linear shape functions are used to approximate the surface, the flow
potential values, and the normal component of its surface gradient. The resulting
Boundary Element Method (BEM, see \cite{brebbia}) consists in collocating a Boundary
Intergal Equation (BIE) in correpondence with each node of the numerical grid, and computing
the integrals appearing in such equation by means of the described iso-parametric formulation.
The linear algebraic equations obtained from such discretization method are combined
with the Ordinary Differential Equations (ODE) derived from the Finite Element spatial 
discretization of the fully nonlinear free surface boundary conditions in Eqs.~\eqref{eq:fsKinematicBeck} and
\eqref{eq:fsDynamicBeck}. The spatial discretization described is carried out making use of the
deal.II open source library for C++ implementation of finite element
discretizations (\cite{BangerthHartmannKanschat2007,BangerthHeisterHeltai-2016-b}). 
To enforce a strong coupling between the fluid and structural problem,
the aforementioned system of Differential Algebraic Equations (DAE) is complemented by the equations
of the rigid hull dynamics (Eqs.~\eqref{eqHullLinCons}-\eqref{eqHullAngCons} and \eqref{eqQuaterEvo}).
The DAE system solution is time integrated by means of an arbitrary
order and arbitrary time step 
implicit Backward Difference Formula (BDF) scheme implemented in the IDA package of the
open source C++ library SUNDIALS (\cite{sundials2005}). The potential flow model described
has been implemented in a stand alone C++ software, the main features of which are
described in~\cite{molaEtAl2013}.

\begin{figure}[htb]
\includegraphics[height=0.28\linewidth]{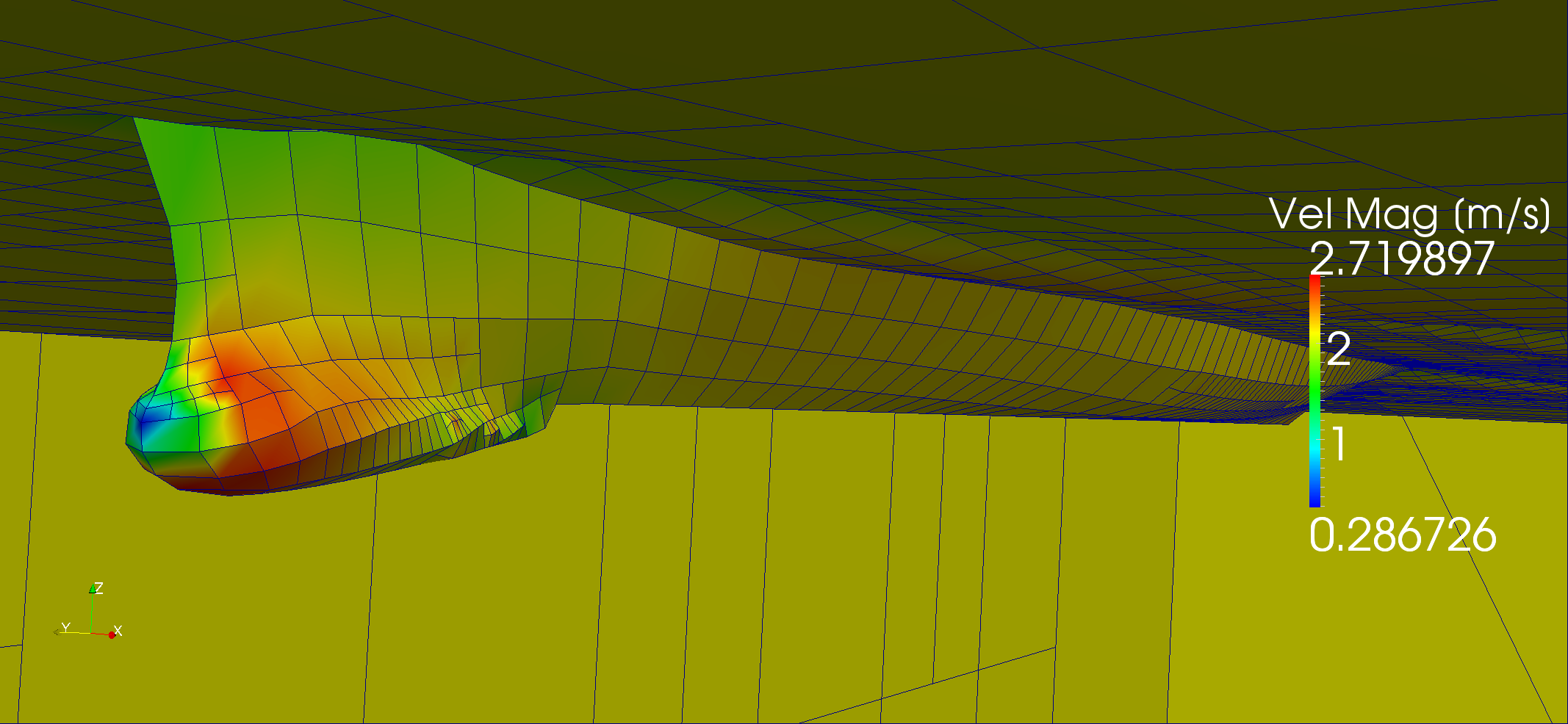}
\includegraphics[height=0.29\linewidth]{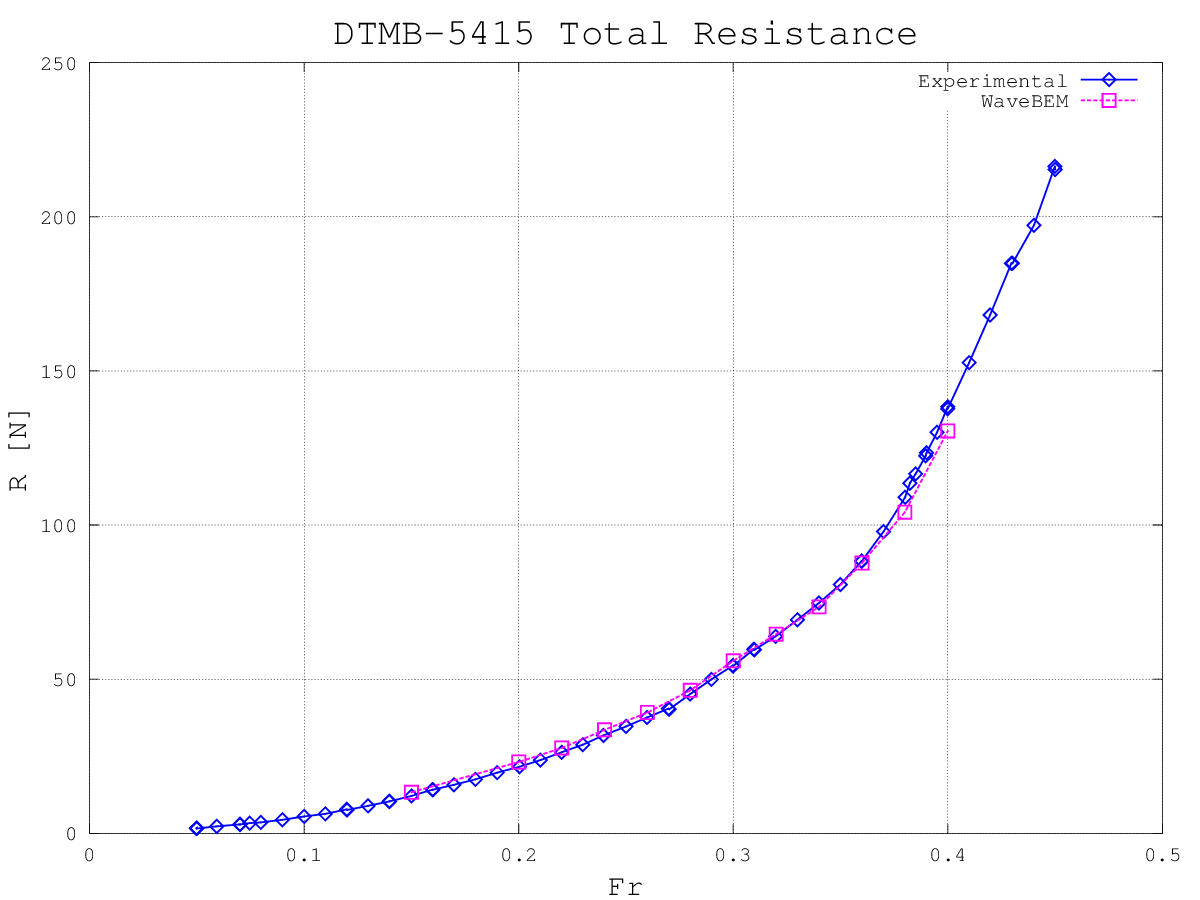}

\caption{On the left, the mesh automatically generated on the surface of one
         of the deformed hulls obtained starting from the
         DTMB 5415 Navy Combatant Hull. On the right, the total resistance
         of the DTMB 5415 hull as a function of the Froude number associated with
         the surge velocity imposed in the simulations. The blue continuous line
         represents the experimental values presented in Olivieri et
         al.~\cite{olivieri2001towing}. The values obtained in this work are
         represented by the dashed magenta line.} \label{fig:waveBemResults}
\end{figure}

The solver is complemented with a mesh module directly interfaced
with CAD data structures~\cite{mola2016ship}. Such feature allows for fully automated mesh
generation once a hull shape is assigned at the start of each simulation by means of
a --- possibly non water-tight --- CAD geometry. Figure~\ref{fig:waveBemResults}, on the
left, displays the mesh generated on the surface of a DTMB 5415 Navy Combatant hull.
At each time step of the simulation, the wave resistance is computed as
\begin{equation}
\label{eqLinForce}
R^w = \int_{\Gamma^b}p\nb\,d\Gamma\cdot\eb_X,
\end{equation}
where $p$ is the pressure value obtained introducing the computed potential in
Eq.~\eqref{eq:bernoulli}. The inviscid fluid dynamic model drag prediction
is then corrected by adding a viscous drag contribution obtained by means of the
ITTC-57 formula~\cite{morrall19701957}. A full assesment of the accuracy of
the high fidelity fluid structure interaction sover described is clearly beyond the scope
of the present work (again, we refer the interested
reader to~\cite{molaEtAl2013,mola2016ship,MolaHeltaiDeSimone2017} for more details).
Yet, in Figure~\ref{fig:waveBemResults}, on the right,
we present a comparison between the computed total resistance curve and the corresponding one
measured by Olivieri et al.~\cite{olivieri2001towing}. As it can be appreciated in
the plot, for all the Froude numbers tested the computed total drag difference
with respect to the measurements is less then 6\%. Given the fact that all the geometries
tested are deformations of the present hull, and that all the velocities imposed fall
in the range reported in the plot, it is reasonable to infer that for
each simulation carried out the accuracy of the high fidelity model prediction will
be similar to that of the results presented.

\section{Parameter space reduction by Active Subspaces}
\label{sec:active}
The active subspaces (AS) approach represents one of the emerging ideas for
dimension reduction in the parameter studies and it is based on the homonymous properties. The concept was
introduced by Paul Constantine in~\cite{constantine2015active}, for example, and employed
in different real world problems. We mention, among others, aerodynamic
shape optimization~\cite{lukaczyk2014active}, the parameter reduction
for the HyShot II scramjet model~\cite{constantine2015exploiting},
active subspaces for integrated hydrologic
model~\cite{jefferson2016reprint}, and in combination with
POD-Galerkin method in cardiovascular problems~\cite{tezzele2017combined}. 

\begin{figure}[h!]
\centering
\subfloat[Original output function]{\includegraphics[width=0.23\textwidth]{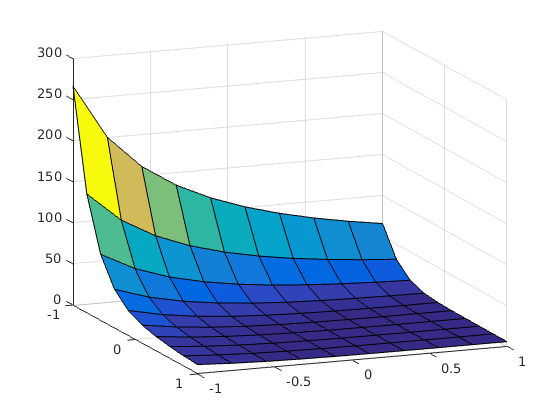}\label{subfig:output_funct}}\hspace{.1cm}
\subfloat[First rotation]{\includegraphics[width=0.23\textwidth]{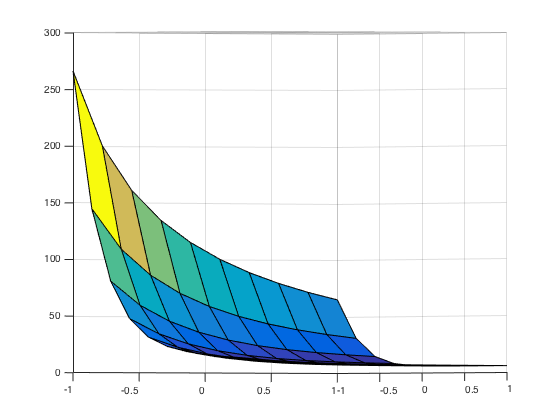}\label{subfig:output_funct_rot1}}\hspace{.1cm}
\subfloat[Second rotation]{\includegraphics[width=0.23\textwidth]{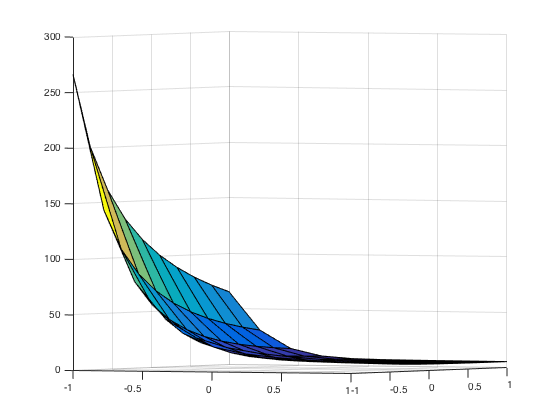}\label{subfig:output_funct_rot2}}\hspace{.1cm}
\subfloat[Output function with respect to the active variable]{\includegraphics[width=0.23\textwidth]{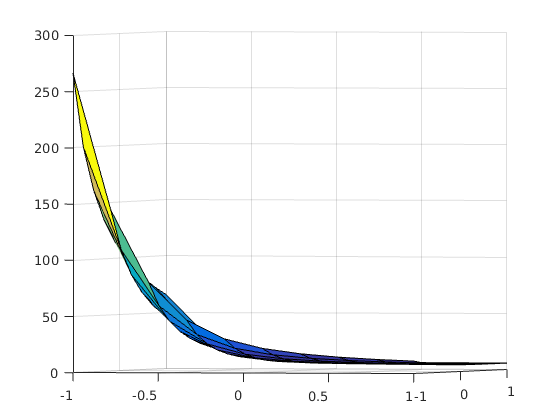}\label{subfig:output_funct_active}}
\caption{Example of a bivariate output function
  \protect\subref{subfig:output_funct}, intermediate rotations of the
  domain \protect\subref{subfig:output_funct_rot1} and \protect\subref{subfig:output_funct_rot2}, and the final
  state \protect\subref{subfig:output_funct_active}, where we can see
  the variation of the function along the active variable.}
\label{fig:active_subspaces_idea}
\end{figure}

A characteristic of the active subspaces is that they identify a set of important
directions in the space of all inputs, instead of identifying a
subset of the inputs as important. If the simulation output
does not change as the inputs move along a particular direction, then
we can safely ignore that direction in the parameter study. In Figure~\ref{fig:active_subspaces_idea} it is
possible to capture the main idea behind the active subspaces
approach: we try to rotate the inputs domain in such a way lower
dimension behavior of the output function is revealed. When an active
subspace is identified for the problem of interest, then it is
possible to perform different parameter studies such as response
surfaces~\cite{box1987empirical}, integration, optimization and statistical inversion~\cite{kaipio2006statistical}.

There are some main ingredients in order to employ active subspaces. The
first is a scalar function $f: \mathbb{R}^m \rightarrow \mathbb{R}$
smooth enough depending on the inputs $\boldsymbol{\mu}$ that represents the
quantity of interest. Moreover we need the gradients of this map with
respect to the inputs
$\nabla_{\boldsymbol{\mu}} f(\boldsymbol{\mu})$ or an approximation of them, a sampling
density $\rho$, and a gap between eigenvalues of the covariance matrix
associated to the gradients. That is a symmetric, positive
semidefinite matrix whose elements are the average products of partial
derivatives of the simulations' input/output map, that reads
\begin{equation}
\label{eq:covariance}
\boldsymbol{\Sigma} = \mathbb{E}\, [\nabla_{\boldsymbol{\mu}} f \, \nabla_{\boldsymbol{\mu}} f
^T] = \int_{\mathbb{D}} (\nabla_{\boldsymbol{\mu}} f) ( \nabla_{\boldsymbol{\mu}} f )^T
\rho \, d \boldsymbol{\mu} ,
\end{equation}
where $\mathbb{E}$ is the expected value, $\rho : \mathbb{R}^m
\rightarrow \mathbb{R}^+$ a probability density
function --- usually a uniform one ---, 
and $\boldsymbol{\Sigma}$ the so-called uncentered covariance matrix of the gradients of $f$ (for a
more deep understanding of these operators see for example~\cite{devore2015probability}).
Usually a Monte Carlo method (see~\cite{metropolis1949monte}) is used in order to approximate the eigenpairs of this
matrix (see \cite{constantine2015computing}) as follows
\begin{equation}
\label{eq:covariance_approx}
\boldsymbol{\Sigma} \approx \frac{1}{N_{\text{train}}^{\text{AS}}} \sum_{i=1}^{N_{\text{train}}^{\text{AS}}} \nabla_{\boldsymbol{\mu}} f_i \,
\nabla_{\boldsymbol{\mu}} f^T_i ,
\end{equation}
where we draw $N_{\text{train}}^{\text{AS}}$ samples $\boldsymbol{\mu}^{(i)}$ independently from the
measure $\rho$ and where $\nabla_{\boldsymbol{\mu}} f_i = \nabla_{\boldsymbol{\mu}}
f(\boldsymbol{\mu}^{(i)})$. This matrix
is symmetric and positive semidefinite, so it admits a real eigenvalue
decomposition
\begin{equation}
\label{eq:decomposition}
\boldsymbol{\Sigma} = \mathbf{W} \mathbf{\Lambda} \mathbf{W}^T ,
\end{equation}
where $\mathbf{W}$ is a $m \times m$ column matrix of eigenvectors,
and $\mathbf{\Lambda}$ is a diagonal matrix of eigenvalues. Then we
order the eigenpairs in descending order. We will select the first $M$
eigenvectors to form a reduced-order basis. This is where we reduce
the dimensionality of the design problem. We will attempt to describe
the behaviour of the objective function by projecting the full-space
design variables onto this active subspace. On average, perturbations
in the first set of coordinates change $f$ more than perturbations in
the second set of coordinates. Low eigenvalues suggest
that the corresponding vector is in the nullspace of the covariance
matrix, and to form an approximation we can discard these vectors. We
select the basis as follow
\begin{equation}
\label{eq:active_decomposition}
\mathbf{\Lambda} =   \begin{bmatrix} \mathbf{\Lambda}_1 & \\
                                     &
                                     \mathbf{\Lambda}_2\end{bmatrix},
\qquad
\mathbf{W} = \left [ \mathbf{W}_1 \quad \mathbf{W}_2 \right ],
\end{equation}
where $\mathbf{\Lambda}_1 = diag(\lambda_1, \dots, \lambda_M)$ with $M<m$, and
$\mathbf{W}_1$ contains the first $M$ eigenvectors. The active subspace
is the span of the first few eigenvectors of the covariance matrix. We
define the active variables to be the linear combinations of the input
parameters with weights from the important eigenvectors. In particular
we define the active subspace to be the range of $\mathbf{W}_1$. The inactive
subspace is the range of the remaining eigenvectors in
$\mathbf{W}_2$. With the basis identified, we can map forward to the
active subspace. So $\boldsymbol{\mu}_M$ is the active variable and
$\boldsymbol{\nu}$ the inactive one, respectively:
\begin{equation}
\label{eq:active_var}
\boldsymbol{\mu}_M = \mathbf{W}_1^T\boldsymbol{\mu} \in \mathbb{R}^M, \qquad
\boldsymbol{\nu} = \mathbf{W}_2^T \boldsymbol{\mu} \in \mathbb{R}^{m-M} .
\end{equation}
In particular any point in the parameter space $\boldsymbol{\mu} \in
\mathbb{R}^m$ can be expressed in terms of $\boldsymbol{\mu}_M$ and
$\boldsymbol{\nu}$: 
\[
\boldsymbol{\mu} = \mathbf{W}\mathbf{W}^T\boldsymbol{\mu} =
\mathbf{W}_1\mathbf{W}_1^T\boldsymbol{\mu} +
\mathbf{W}_2\mathbf{W}_2^T\boldsymbol{\mu} = \mathbf{W}_1 \boldsymbol{\mu}_M +
\mathbf{W}_2 \boldsymbol{\nu} .
\]
Having this decomposition in mind we can rewrite $f$
\[ 
f (\boldsymbol{\mu}) =  f (\mathbf{W}_1 \boldsymbol{\mu}_M + \mathbf{W}_2
 \boldsymbol{\nu}) ,
\]
and construct a surrogate model $g$ discarding the inactive variables
\[
f (\boldsymbol{\mu}) \approx g (\mathbf{W}_1^T \boldsymbol{\mu}) = g(\boldsymbol{\mu}_M).
\]

\TODO{The surrogate model $g$ in this work is a response
surface. We exploit the decreased number of parameters to fit a lower
dimensional response surface, fighting the curse of dimensionality
that affects this approximation procedure. The advantage of this
approach is that more models are feasible, such as for example radial
basis functions interpolation, higher degree polynomials, or
regressions techniques. In particular we use a polynomial response
surface. For different type of surrogate model that exploit a shared
active subspaces in naval engineering refer to~\cite{tezzele2018shared}.}

We underline that the size of the eigenvalue problem is the limiting factor. We need to
compute eigenvalue decompositions with $m \times m$ matrices, where
$m$ is the dimension of the simulation, that is the number of inputs. 

Active subspaces can be seen in the more general context of ridge
approximation (see \cite{pinkus2015ridge,keiper2015analysis}). In particular it can be proved that, under certain
conditions, the active subspace is nearly stationary and it is a good
starting point in optimal ridge approximation as shown in \cite{constantine2016near,hokanson2017data}.

\section{Numerical results}
\label{sec:results}
In this section we present the results of the complete pipeline, 
presented in the previous sections and in Figure~\ref{fig:scheme}, 
applied to the DTMB 5415 hull.

The mesh is discretized with quadrilateral cells. The BEM uses
bi-linear quadrilateral elements. This results in roughly 4000 degrees
of freedom for each simulation realized. The high
fidelity solver described in section~\ref{sec:bem} is implemented
in WaveBEM~\cite{MolaHeltaiDeSimone2017} using the deal.II
library~\cite{BangerthHeisterHeltai-2016-b}. 

\begin{figure}[h!]
\centering
\subfloat[Original wave resistance]{\includegraphics[width=0.45\textwidth]{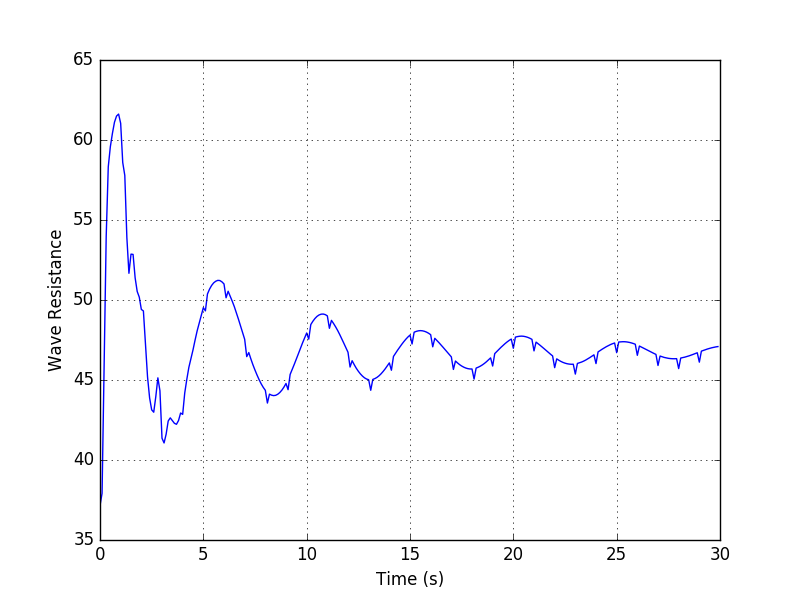}\label{subfig:wave_res}}\hspace{.5cm}
\subfloat[Fitted wave resistance]{\includegraphics[width=0.45\textwidth]{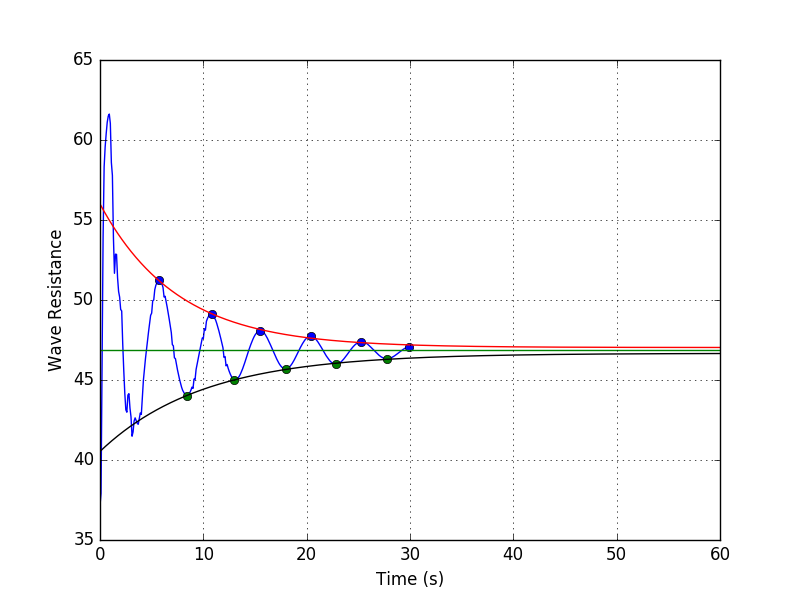}\label{subfig:wave_res_fitted}}
\caption{Plot \protect\subref{subfig:wave_res} shows the original wave
  resistance simulated for 30 seconds. Then plot \protect\subref{subfig:wave_res_fitted}
  depicts, after a filter has been applied, the exponential fitting of the maximums and minimums and the average at regime for 60 seconds.}
\label{fig:original_wave}
\end{figure}

The solver, after reading the CAD file of the deformed geometry,
simulates the behaviour of the hull for 30 seconds. To further speedup the computations, the
total resistance computed as in Eq.~\eqref{eqLinForce} is
extrapolated to obtain the total resistance at the final, steady
state regime, with an error in the order of 0.1\%. In
Figure~\ref{fig:original_wave} we can see the original simulation of
the total resistance for the first 30 seconds and then the
extrapolation we have done after a proper filtering of the data. We
fit the maximums using the following function: $a e^{-b x} +
c$. For the minimums we use $-a e^{-b x} + c$. Then we set the approximated
wave resistance to the average of the two at infinity.

Let us recall that the parameter space is a $m = 8$ dimensional
space. The parameters are showed in
Table~\ref{tab:design_space}. We remark that the first six are
geometrical parameters that produce the deformation of the original
domain, while the last two are structural and physical parameters ---
the displacement and the velocity of the hull ---. The PyGeM open
source package is used to perform the free form
deformation~\cite{pygem}. 

We create a dataset with 130 different couples of input/output data. We split the dataset in two,
creating a train dataset with 80\% of the data, and a test dataset
with the remaining~20\%. That means that
$N_{\text{train}}^{\text{AS}} = 104$ in
Eq.~\eqref{eq:covariance_approx}. Even though it may be challenging to
explore a 8 dimensional space, as reported
in~\cite{constantine2015active}, heuristics suggest that this choice
of $N_{\text{train}}^{\text{AS}}$ is enough for the purposes of the
active subspaces identification described in section~\ref{sec:active}.

In order to construct the uncentered covariance matrix $\boldsymbol{\Sigma}$,
defined in Eq.~\eqref{eq:covariance}, we use a Monte Carlo method as
shown in Eq.~\eqref{eq:covariance_approx}, employing the software in~\cite{paul-contantine-16}. Since we have only pairs of input/output
data we need to approximate the gradients of the total wave resistance
with respect to the parameters, that is $\nabla_{\boldsymbol{\mu}}
f$. We use a local linear model that approximates the gradients with the best linear
approximation using 14 nearest neighbors. After constructing the
matrix we calculate its real eigenvalue decomposition. Recalling
section~\ref{sec:active},since $m = 8$, we have that $\boldsymbol{\Sigma} \in
\mathbb{M} (8, \mathbb{R})$.

\begin{figure}[h]
\centering
\subfloat[Eigenvalues estimates.]{\includegraphics[width=0.43\textwidth]{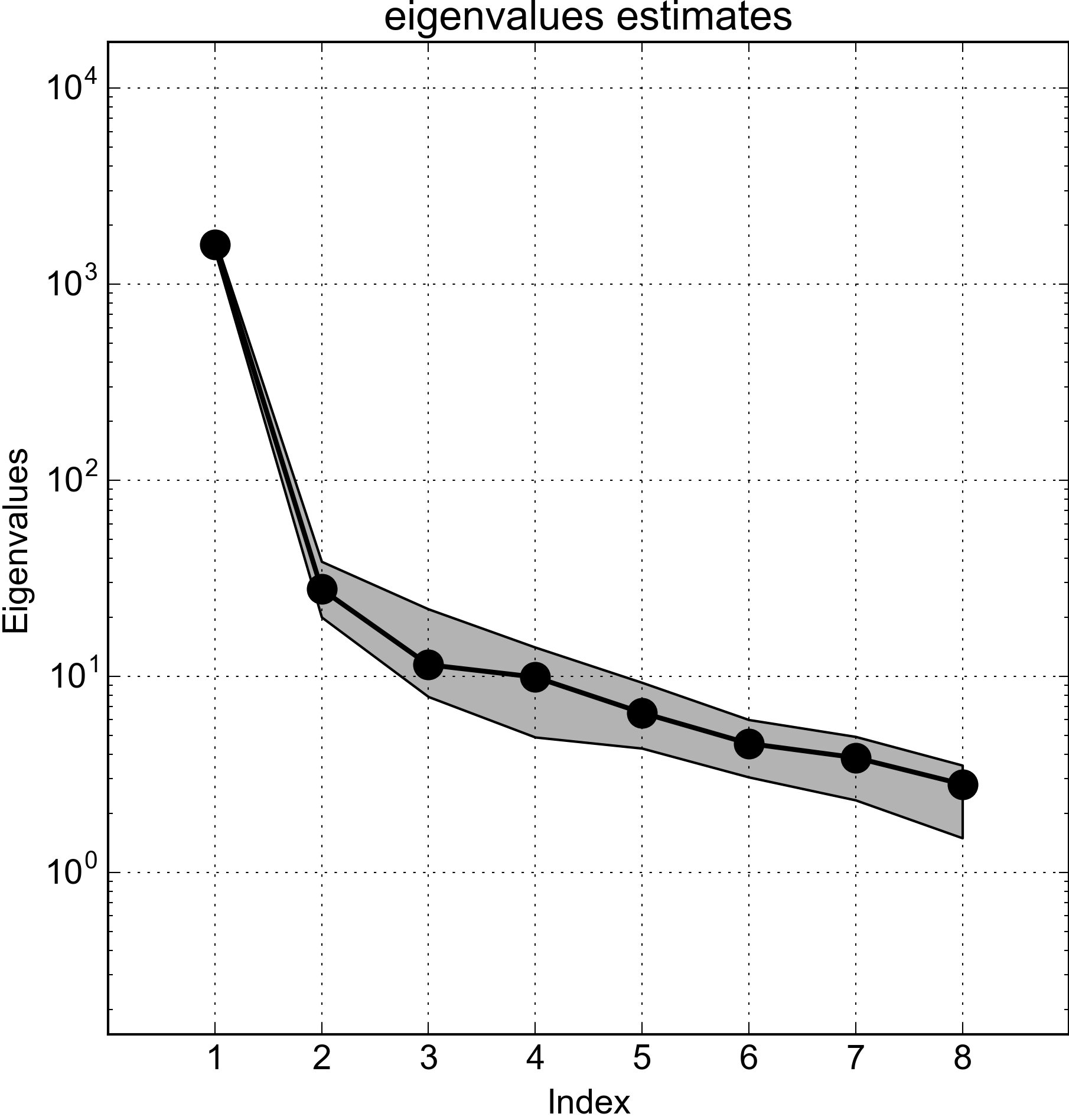}\label{subfig:eigenvalues_hull}}\hspace{.5cm}
\subfloat[First eigenvector estimate.]{\includegraphics[width=0.47\textwidth]{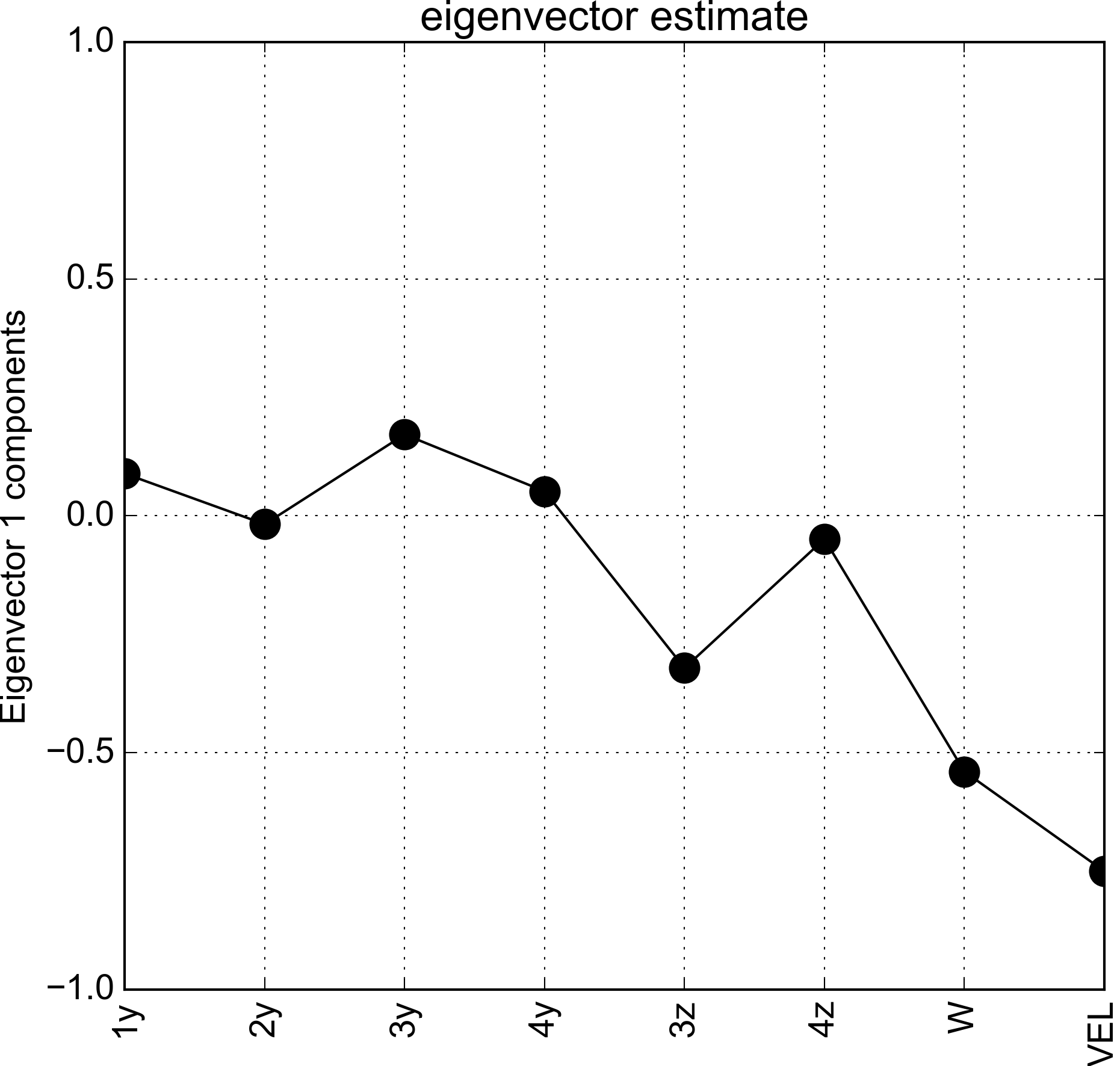}\label{subfig:eigenvector_hull}}
\caption{Plot \protect\subref{subfig:eigenvalues_hull}
  shows the eigenvalue estimates in block circles with the bootstrap
  intervals (grey region). The
  order-of-magnitude gaps between the eigenvalues suggest confidence
  in the dominance of the active
  subspace. Plot \protect\subref{subfig:eigenvector_hull} shows the
  components of the eigenvector correspondent to the greatest 
  eigenvalue.}
\label{fig:eigen_hull}
\end{figure}

In Figure~\ref{subfig:eigenvalues_hull} we see the eigenvalues of $\boldsymbol{\Sigma}$
and the bootstrap intervals. Bootstrapping
is the practice of estimating properties of a quantity (such as its
variance) by measuring those properties when sampling from an
approximating distribution. It relies on random sampling with
replacement. The bootstrap intervals in grey are computed using 1000
bootstrap replicates randomly extracted from the original gradient
samples. We underline the presence of a major gap between the first and the
second eigenvalue and a minor one between the second and the
third. This suggests the existence of a one dimensional subspace and
possibly the presence of a two dimensional one. To better investigate
the first case, in Figure~\ref{subfig:eigenvector_hull} we present the
components of the eigenvector with index 1 that corresponds to
the greatest eigenvalue, that is the matrix --- in this case it is a
vector --- $\mathbf{W}_1 \in
\mathbb{R}^8$ of Eq.~\eqref{eq:active_decomposition}. Since they are the weights of the linear
combination used to construct the active direction they provide useful
informations. We can see that the major contributions are due to the
velocity, the weight, and the depth of the hull. We underline that a
weight that is almost zero means that the output function, on average, is almost
flat along the direction identified by the corresponding
parameter. This is a very useful information for a designer because in
such a way he can deform the shape along prescribed directions without
affecting the total wave resistance, allowing for example to store
more goods inside the hull preserving the performances.

\begin{figure}[h]
\centering
\subfloat[One active variable.]{\includegraphics[width=0.44\textwidth]{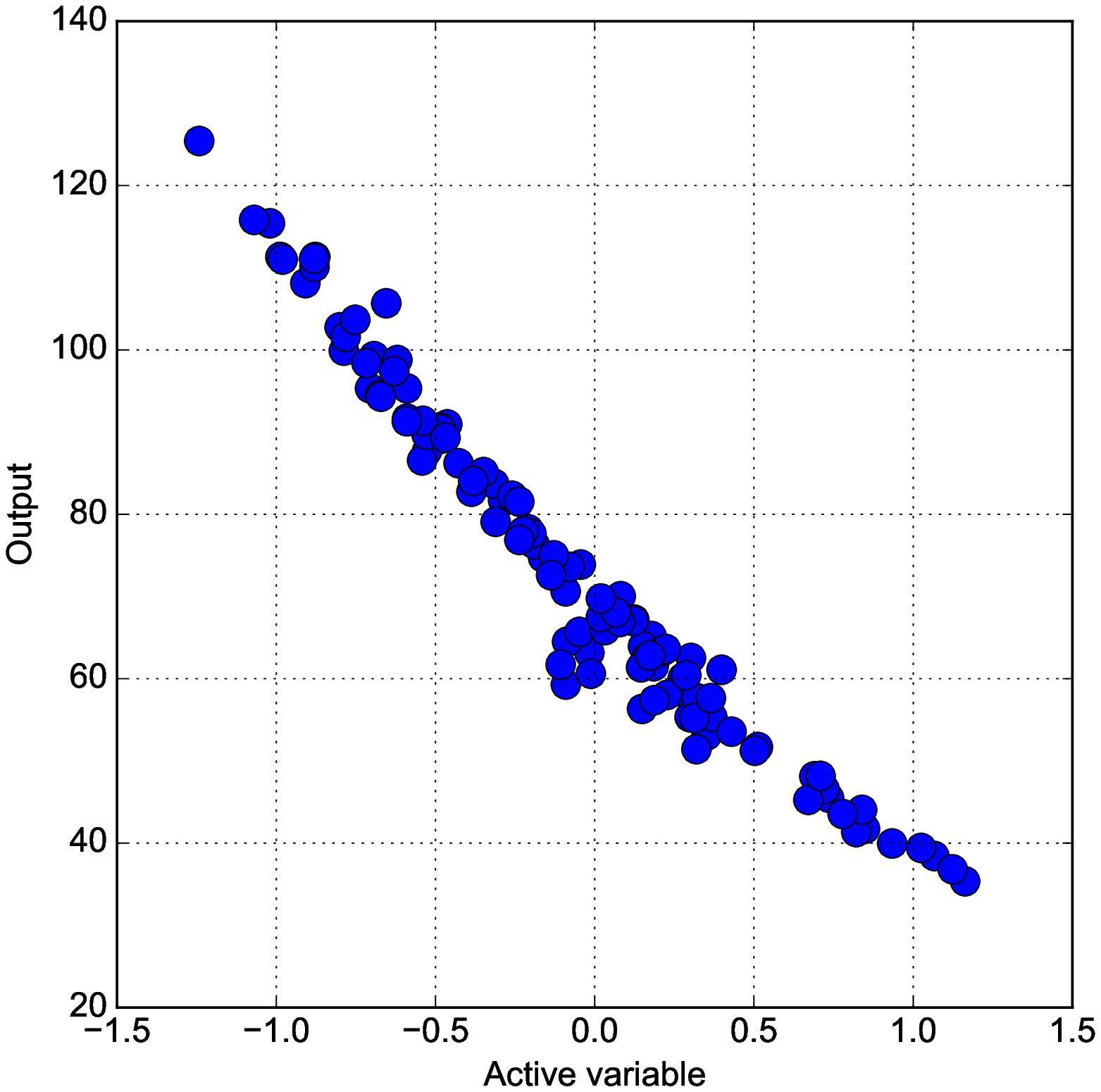}\label{subfig:ssp1_hull}}\hspace{.5cm}
\subfloat[Two active variables.]{\includegraphics[width=0.48\textwidth]{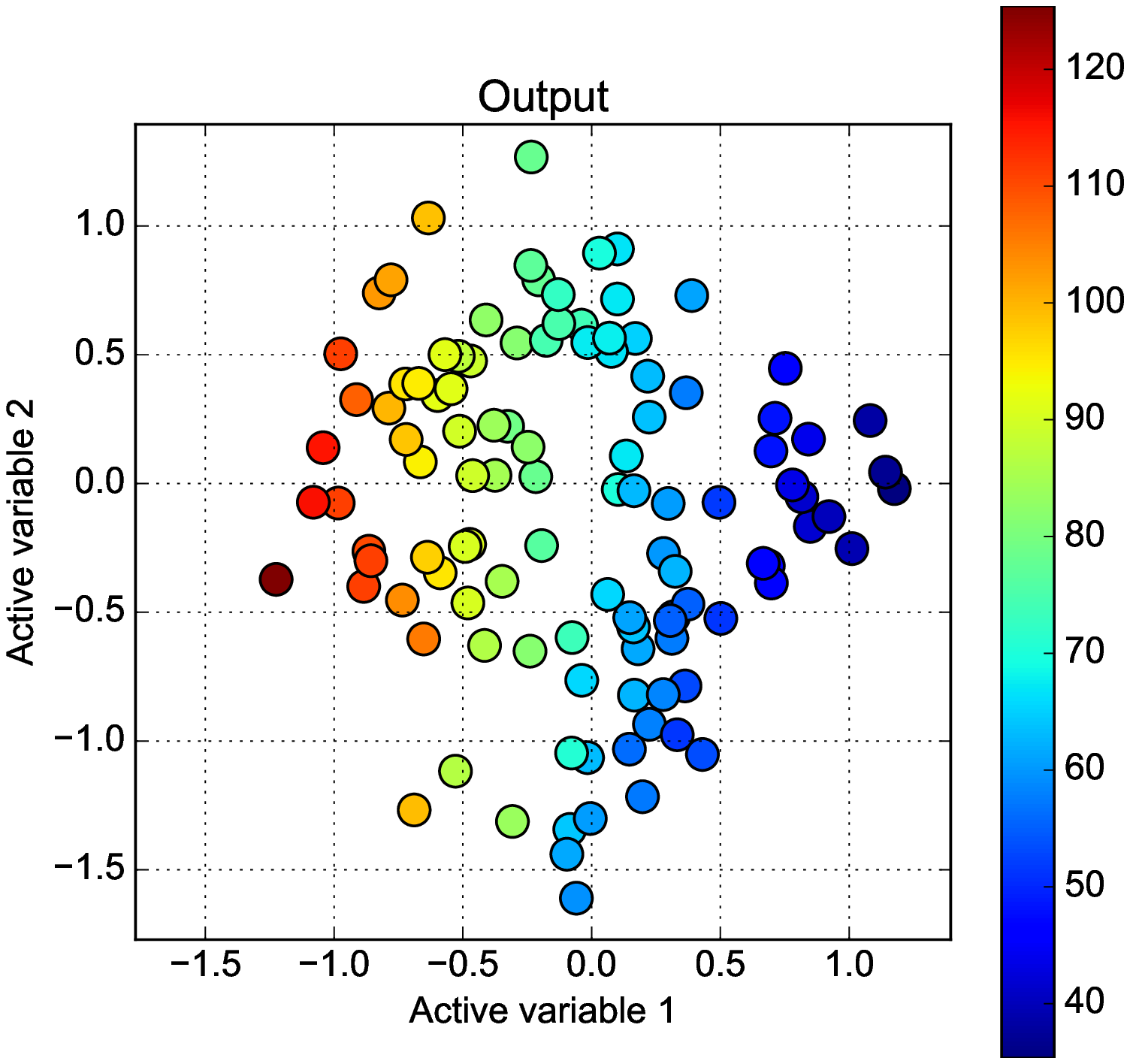}\label{subfig:ssp2_hull}}
\caption{Sufficient summary plots for \protect\subref{subfig:ssp1_hull} one and \protect\subref{subfig:ssp2_hull} two active variables using the training dataset.}
\label{fig:ssp_hull}
\end{figure}

In Figure~\ref{fig:ssp_hull} we explore the training dataset,
plotting the sufficient summary plot (see \cite{cook2009regression}) for one and two active
variables. Sufficient summary plots are powerful visualization tools
to identifying low-dimensional structure in a quantity that depends on
several input variables. A scatter plot that contains all available
regression information is called a sufficient summary plot. Recalling
Eq.~\eqref{eq:active_var}, Figure~\ref{fig:ssp_hull} shows
$f(\boldsymbol{\mu})$ against $\boldsymbol{\mu}_M = \mathbf{W}_1^T \boldsymbol{\mu}$, where
$\mathbf{W}_1$ contains the first one and two eigenvectors
respectively. In particular each point represents the value of the
output function for a particular choice of the parameters, mapped in
the active subspace. The two plots confirm the presence of an active subspace of
dimension one and two. The latter seems to capture the output function
in a much finer way, but as we are going to show the gain in terms of
error committed is not so big.

\begin{figure}[htb]
\centering
\subfloat[Surrogate model error with respect to the active subspace
dimension and the order of the response
surface.]{\includegraphics[width=0.56\textwidth]{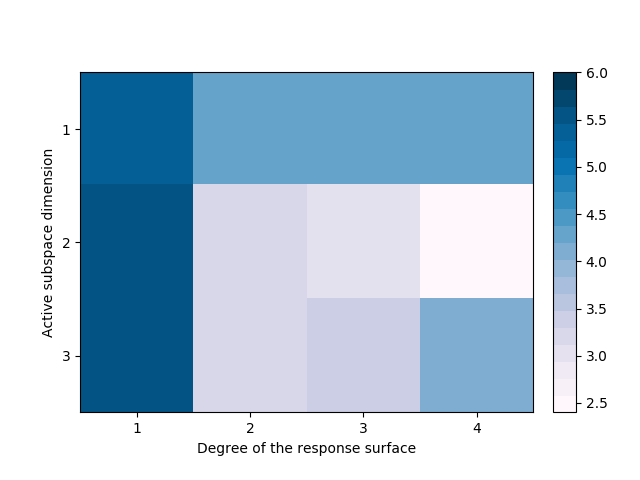}\label{subfig:heat_map}}\hspace{.5cm}
\subfloat[Observations and the corresponding predictions using a
polynomial response surface of order three and dimension one.]{\includegraphics[width=0.385\textwidth]{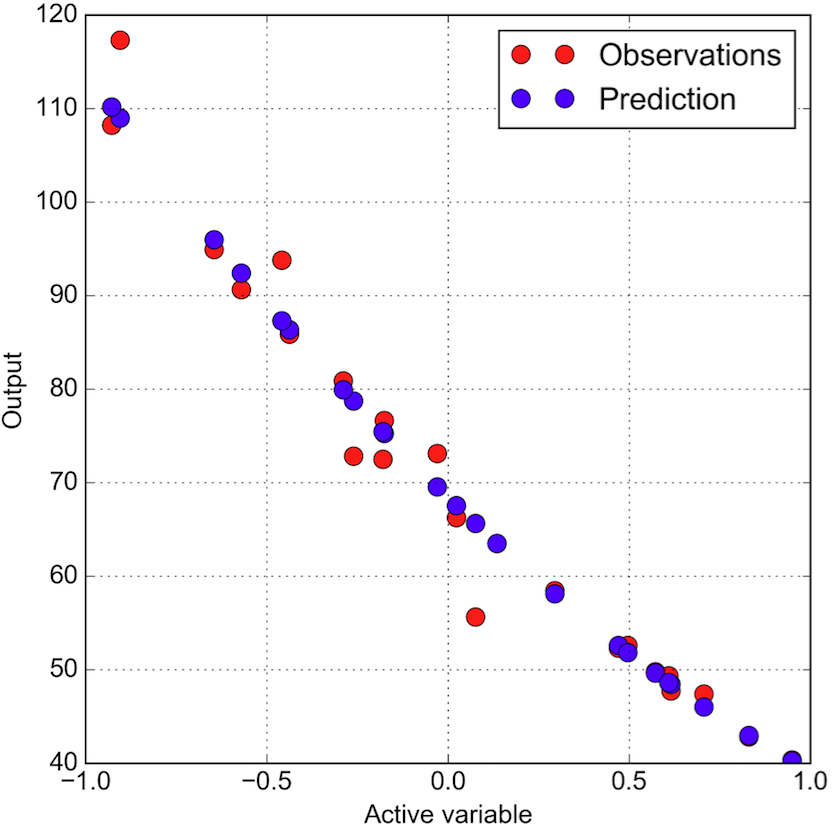}\label{subfig:predictions_deg3_dim1}}
\caption{Plot \protect\subref{subfig:heat_map}
  shows the relative error on the test dataset with respect to the active subspace dimension and the order of the response
surface; plot \protect\subref{subfig:predictions_deg3_dim1} shows the observations of the test dataset and the corresponding
  predictions using a polynomial response surface of order three.}
\label{fig:hull_errors_deg}
\end{figure}

We can compare the decay of the eigenvalues with the decay of
surrogate model error on the test dataset shown in Figure~\ref{subfig:heat_map}. We project the data
onto active subspaces of varying dimension, and construct a surrogate
model with a least-squares-fit, global, multivariate polynomial
approximation of different orders. Then we calculate the root-mean-square
error of the test data against the surrogate. This error is
scaled with respect to the range of the function evaluations, making
it a relative error. We repeat this procedure 20 times constructing
every time the uncentered covariance matrix of Eq.~\eqref{eq:covariance_approx},
since a Monte Carlo approximation is involved. Finally we take the
average of the errors computed. Because we have a large amount of training data,
we can expect the surrogate model constructed in a low dimension to be
accurate if the data collapses into a manifold. Thus the test
error is an indication of how well the active subspace has collapsed
the data. In Figure~\ref{subfig:heat_map} are depicted the errors with
respect to the active subspace dimension and the order of the response
surface. The subspace dimension varies from 1 to 3, while the order of
the response surface from 1 to 4. The minimum error is achieved using
a two dimensional active variable and a response surface of degree~4
and it is $\approx 2.5\%$. Further investigations show that increasing
the dimension of the active variable does not decrease significantly
the error committed while the time to construct the corresponding response surface
increases. This is confirmed by the marginal gains in the decay of the
eigenvalues for active variables of dimension greater then three as
shown in Figure~\ref{subfig:eigenvalues_hull}. 
We can affirm that the active
subspace of dimension one is sufficient to model the wave resistance
of the DTMB 5415 if we can afford an error of approximately 4.5\%. In particular in
Figure~\ref{subfig:predictions_deg3_dim1} we can see the predictions
made with the surrogate model of dimension one and the actual
observations. Otherwise, we can achieve a $\approx 2.5\%$ error if we
take advantage of two active dimensions and a response surface of order
four, preserving a fast evaluation of the surrogate model.

We want to stress the fact that the result is remarkable if we
consider the heterogeneous nature of all the parameters involved. In the
case of only geometrical parameters one can easily expect such a
behaviour but considering also physical and structural ones make the
result not straightforward at all. Moreover the evaluation of the response
surface takes less than one second compared to the 12 hours of a full
simulation per single set of parameters on the same computing
machine. This opens new potential approaches to optimization problems.

\section{Conclusions and perspectives}
\label{sec:the_end}
In this work we presented a numerical framework for the reduction of
the parameter space and the construction of an optimized response surface to
calculate the total wave resistance of the DTMB 5415 advancing in calm
water. We integrate heterogeneous parameters in order to have insights
on the more important parameters. The reduction both in terms of cost
and time, remaining below the 4.2\% of error on new unprocessed data,
is very remarkable and promising as a new design interpreted tool. The methodological and computational pipeline is also easily extensible to different hulls
and/or different parameters. This allows a fast preprocessing and a very good
starting point to minimize the quantities of interest in the field of optimal shape design.

This work is directed in the development of reduced order methods
(ROMs) and efficient parametric studies. Among others we would like to
cite~\cite{chinesta2016model,hesthaven2016certified,salmoiraghi2016advances,salmoiraghi2017}
for a comprehensive overview on ROM and geometrical
deformation. Future developments involve the application of the POD
after the reduction of the parameter space through the active
subspaces approach.

\section*{Acknowledgements}
This work was partially supported by the project OpenViewSHIP,
``Sviluppo di un ecosistema computazionale per la progettazione
idrodinamica del sistema elica-carena'' and ``Underwater Blue Efficiency'', supported by Regione FVG -
PAR FSC 2007-2013, Fondo per lo Sviluppo e la Coesione, by the
project ``TRIM – Tecnologia e Ricerca Industriale per la Mobilit\`a
Marina'', CTN01-00176-163601, supported by MIUR, the italian
Ministry of Education, University and Research, by the INDAM-GNCS 2017
project ``Advanced numerical methods combined with computational
reduction techniques for parameterised PDEs and applications'' and by
European Union Funding for Research and Innovation --- Horizon 2020
Program --- in the framework of European
Research Council Executive Agency: H2020 ERC CoG 2015 AROMA-CFD
project 681447 ``Advanced Reduced Order Methods with Applications in
Computational Fluid Dynamics'' P.I. Gianluigi Rozza.



\end{document}